\documentclass[12pt]{article}

\usepackage{amsmath, epsfig, cite}
\usepackage{amssymb}
\usepackage{amsfonts}
\usepackage{latexsym}
\usepackage{graphicx}
\usepackage{enumerate}

\newtheorem{thm}{Theorem}[section]

\newtheorem{cor}[thm]{Corollary}

\newtheorem{lem}[thm]{Lemma}

\newtheorem{problem}[thm]{Problem}

\newcommand{\qed}{{\hfill\rule{4pt}{7pt}}}
\def\pf{\noindent {\it Proof.} }

\numberwithin{equation}{section}

\makeatletter \@addtoreset{equation}{section} \makeatother

\title {\bf Rainbow $k$-connectivity of\\ random bipartite graphs\footnote{Supported by NSFC No. 11071130 and the ``973'' project.}}

\author{
{\small Xiaolin Chen, Xueliang Li, Huishu Lian}\\
{\small Center for Combinatorics and LPMC-TJKLC}\\
{\small Nankai University, Tianjin 300071, P.R. China}\\
{\small E-mail: chxlnk@163.com; lxl@nankai.edu.cn; lhs6803@126.com}
   }
\date{}
\begin{document}

\maketitle

\begin{abstract}
A path in an edge-colored graph $G$ is called a \emph{rainbow path}
if no two edges of the path are colored the same. The minimum number
of colors required to color the edges of $G$ such that every pair of
vertices are connected by at least $k$ internally vertex-disjoint
rainbow paths is called the \emph{rainbow $k$-connectivity} of the
graph $G$, denoted by $rc_k(G)$. For the random graph $G(n,p)$, He
and Liang got a sharp threshold function for the property
$rc_k(G(n,p))\leq d$. In this paper, we extend this result to
the case of random bipartite graph $G(m,n,p)$.\\

\noindent\textbf{Keywords:} rainbow $k$-connectivity, sharp threshold function,
random bipartite graph\\

\noindent\textbf{AMS Subject Classification Numbers:} 05C15, 05C40,
05C80
\end{abstract}

\section{Introduction}

In this paper, unless otherwise stated, all graphs are finite,
simple and undirected. For basic terminology and notation in graph
theory, see\cite{JU}. Connectivity is one of the basic concepts of
graph theory. Recently, the concepts of rainbow connectivity (or
rainbow connection) and rainbow $k$-connectivity are introduced by
Chartrand et. al. in \cite{GGKP1} and \cite{GGKP2} as a
strengthening of the canonical connectivity concept. Given an
edge-colored graph $G$, we call a path a \emph{rainbow path} if no
two edges of the path are colored the same. We call the graph $G$
\emph{rainbow connected} if every pair of vertices are connected by
at least one rainbow path. The minimum number of colors required to
make $G$ rainbow connected is called the \emph{rainbow
connectivity}, denoted by $rc(G)$. In general, for an integer $k\geq
1$, a graph $G$ is called \emph{rainbow $k$-connected} if every pair
of vertices of $G$ are connected by at least $k$ internally
vertex-disjoint rainbow paths. The minimum number of colors required
to make $G$ rainbow $k$-connected is called the \emph{rainbow
$k$-connectivity}, denoted by $rc_k(G)$.

In addition to regarding it as a natural combinatorial concept,
rainbow connectivity also has interesting applications in
transferring information of high security and networking
\cite{GGKP2}, \cite{SEAR} and \cite{A}. The following motivation
comes from \cite{SEAR}: Suppose we wish to route messages between
any two vertices in a cellular network and require that each link on
the route between the vertices is assigned with a distinct channel.
We clearly wish to minimize the number of distinct channels. The
minimum number is exactly the rainbow connectivity of the underlying
graph. The subject has since attracted considerable interest. A
great number of results about the rainbow connectivity have been
obtained by the researchers. Recently, Li and Sun published a book
\cite{LiSun} and Li et. al. wrote a survey \cite{LSS} on the current
status of rainbow connectivity. We refer them to the reader for
details.

We will study the rainbow $k$-connectivity in the random graph
setting \cite{B}. Some results have been obtained in the
Erd\H{o}s-R\'{e}nyi random graph model $G(n,p)$, which is a graph
with $n$ vertices where each of the ${n \choose 2}$ potential edges
appears with probability $p$, independently. Random bipartite graph
model is a general model for complex networks, thus in this paper,
we will extend the results to the random bipartite graph $G(m,n,p)$
with bipartition $(U,V)$, where $|U|=m,|V|=n$ and for each $u\in U$
and $v\in V$ the edge $uv$ appears with probability $p$,
independently. We say that an event $E=E(n)$ happens \emph{almost
surely} (or a.s. for brevity) if $\lim_{n\rightarrow
\infty}\text{Pr}[E(n)]=1$. For a graph property $\mathcal{P}$, a
function $p^*(n)$ is called a \emph{threshold function} of
$\mathcal{P}$ if \footnote[2]{We use the following standard
asymptotic notations: as $n\rightarrow \infty$, $f(n)=o(g(n))$ means
that $f(n)/g(n)\rightarrow 0$; $f(n)=\omega(g(n))$ means that
$f(n)/g(n)\rightarrow \infty$; $f(n)=O(g(n))$ means that there
exists a constant $C$ such that $|f(n)|\leq C g(n)$;
$f(n)=\Omega(g(n))$ means that there exists a constant $c>0$ such
that $f(n)\geq c g(n)$;}
\begin{enumerate}[\indent $\bullet$]
\item for every $p(n)=o(p^*(n))$, $G(n,p(n))$ almost surely does not satisfy $\mathcal{P}$; and
\item for every $p'(n)=\omega(p^*(n))$, $G(n,p'(n))$ almost surely satisfies $\mathcal{P}$.
\end{enumerate}
Furthermore, $p^*(n)$ is called a sharp threshold function of
$\mathcal{P}$ if there are two positive constants $c$ and $C$ such
that
\begin{enumerate}[\indent $\bullet$]
\item for every $p(n)\leq c\cdot p^*(n)$, $G(n,p(n))$ almost surely does not satisfy $\mathcal{P}$; and
\item for every $p'(n)\geq C\cdot p^*(n)$, $G(n,p'(n))$ almost surely satisfies $\mathcal{P}$.
\end{enumerate}

It is well known that all monotone graph properties have a sharp
threshold function \cite{BA} and \cite{EG}. Obviously for every
$k,d$, the property that the rainbow $k$-connectivity is at most $d$
is monotone, and thus has a sharp threshold. Caro et. al.
\cite{YAYZR} proved that $p=\sqrt{\log n/n}$ is a sharp threshold
function for the property $rc(G(n,p))\leq 2$. This was generalized
by He and Liang \cite{HL}, who proved that if $d\geq 2$ and $k\leq
O(\log n)$, then $p=(\log n)^{1/d}/n^{(d-1)/d}$ is a sharp threshold
function for the property $rc_k(G(n,p))\leq d$. Moreover, Fujita et.
al. \cite{FLM} proved that in the random bipartite graph $G(n,n,p)$,
$p=\sqrt{\log n/n}$ is a sharp threshold function for the property
$rc_k(G(n,n,p))\leq 3$. They also posed some open problems, one of
which is stated as follows.
\begin{problem} \cite{FLM}
For $d\geq 2$, determine a sharp threshold function for the property
$rc_k(G)\leq d$, where $G$ is another random graph model.
\end{problem}

In this paper, we consider the random bipartite graph $G(m,n,p)$.
The following results are obtained.
\begin{thm}\label{th1}
Let $d\geq 2$ be a fixed positive integer and $k=k(n)\leq O(\log n)$.\\
\indent If $d$ is odd, then
\begin{equation*}
p=(\log (mn))^{1/d}/(m^{(d-1)/(2d)}n^{(d-1)/(2d)})
\end{equation*}
is a sharp threshold function for the property $rc_k(G(m,n,p))\leq d+1$,
where $m$ and $n$ satisfy that $pn\geq pm \geq (\log n)^4$.\\
\indent If $d$ is even, then
\begin{equation*}
p=(\log n)^{1/d}/(m^{1/2}n^{(d-2)/(2d)})
\end{equation*}
is a sharp threshold function for the property $rc_k(G(m,n,p))\leq
d+1$, where $m$ and $n$ satisfy that there exists a small constant
$\epsilon$ with $0<\epsilon<1$ such that $pn^{1-\epsilon}\geq
pm^{1-\epsilon} \geq (\log n)^4$.
\end{thm}

Then, the following corollary follows immediately.
\begin{cor}
Let $d\geq 2$ be a fixed integer and $k=k(n)\leq O(\log n)$.
Then $p=(\log n)^{1/d}/n^{(d-1)/d}$ is a sharp threshold function
for the property $rc_k(G(n,n,p))\leq d+1$.
\end{cor}

When $d=2$, we get the result of Fujita et. al. in \cite{FLM}.

In the sequel, we will first show Theorem \ref{th1} in Section 2.
Then in Section 3, we will prove a conclusion stated in Section 2,
which plays a key role during our proof of Theorem \ref{th1}.

\section{Threshold of the rainbow $k$-connectivity}

In this section, we establish a sharp threshold function of the
random bipartite graph $G(m,n,p)$ for the property
$rc_k(G(m,n,p))\leq d+1$. We distinguish two parts to prove Theorem
\ref{th1} according to the parity of $d$. For brevity, let
$p_1=(\log (mn))^{1/d}/(m^{(d-1)/(2d)}n^{(d-1)/(2d)})$ and
$p_2=(\log n)^{1/d}/(m^{1/2}n^{(d-2)/(2d)})$. For a fixed $d$, we
always assume that $p_1n\geq p_1m \geq (\log n)^4$ if $d$ is odd and
there exists a $0<\epsilon<1$ such that $p_2n^{1-\epsilon}\geq
p_2m^{1-\epsilon} \geq (\log n)^4$ if $d$ is even. In the sequel, we
fix $\epsilon$. Before our proof, we first recall the following fact
on the diameter of a random bipartite graph.
\begin{thm} $[2, p.272]$ \label{th2}
Suppose that for all $n$,
\begin{equation*}
pn\geq pm\geq (\ln n)^4
\end{equation*}
and that $d$ is a fixed positive integer.\\
\indent If $d$ is odd and
\begin{equation*}
p^dm^{(d-1)/2}n^{(d-1)/2}-\ln (mn)\rightarrow \infty ,
\end{equation*}
or if $d$ is even and
\begin{equation*}
p^dm^{d/2}n^{d/2-1}-2\ln n\rightarrow \infty ,
\end{equation*}
then almost every $G(m,n,p)$ is of diameter at most $d+1$.\\
\indent If $d$ is odd and
\begin{equation*}
p^dm^{(d-1)/2}n^{(d-1)/2}-\ln (mn)\rightarrow -\infty ,
\end{equation*}
or if $d$ is even and
\begin{equation*}
p^dm^{d/2}n^{d/2-1}-2\ln n\rightarrow -\infty ,
\end{equation*}
then almost every $G(m,n,p)$ is of diameter at least $d+2$.
\end{thm}

From the above theorem, we can derive that \\
\indent when $d$ is odd, for every $c<1$ and $p(n)\leq c\cdot(\ln
(mn))^{1/d}/(m^{(d-1)/(2d)}n^{(d-1)/(2d)})$, $G(m,n,p)$ almost
surely does not satisfy the property that diam$(G(m,n,p))\leq d+1$
and for every $C>1$ and $p(n)\geq C\cdot(\ln
(mn))^{1/d}/(m^{(d-1)/2d}n^{(d-1)/2d})$, $G(m,n,p)$ almost surely
satisfies the property that diam$(G(m,n,p))\leq d+1$. Similarly,
when $d$ id even, for every $c<1$ and $p(n)\leq c\cdot(2\log
n)^{1/d}/(m^{1/2}n^{(d-2)/(2d)})$, $G(m,n,p)$ almost surely does not
satisfy the property that diam$(G(m,n,p))\leq d+1$ and for every
$C>1$ and $p(n)\geq C\cdot(2\log n)^{1/d}/(m^{1/2}n^{(d-2)/(2d)})$,
$G(m,n,p)$ almost surely satisfies the property that
diam$(G(m,n,p))\leq d+1$.

We also need the following key conclusion during our proof. Here we
only state it but give its proof next section. Assume that
$c_0\geq1$ is a positive constant and $k=k(n)\leq c_0\log n$. Let
$C_1=2^{10d}\cdot c_0$ and $C_2=2^{10d}\cdot c_0/\epsilon$.
\begin{thm}\label{th3}
\indent If $d$ is odd, then with probability at least $1-n^{-\Omega(1)}$,
the random bipartite graph $G(m,n,C_1p_1)$ satisfies the property:\\
\indent every two distinct vertices of the same partite are
connected by at least $2^{10d}c_0\log n$ internally vertex-disjoint
paths of length exactly $d+1$, and every two distinct vertices of
different partites are connected by at least $2^{10d}c_0\log n$
internally vertex-disjoint paths of length exactly $d$.\\
\indent If $d$ is even, then with probability at least $1-n^{-\Omega(1)}$,
the random bipartite graph $G(m,n,C_2p_2)$ satisfies the property:\\
\indent every two distinct vertices of the same partite are
connected by at least $2^{10d}c_0\log n$ internally vertex-disjoint
paths of length exactly $d$, and every two distinct vertices of
different partites are connected by at least $2^{10d}c_0\log n$
internally vertex-disjoint paths of length exactly $d+1$.
\end{thm}

Now we are ready to give the proof of Theorem \ref{th1}.

\noindent\textbf{Part 1: $d$ is odd.}

We consider the random bipartite $G(m,n,p)$ with $p_1n\geq p_1m \geq
(\log n)^4$. To establish a sharp threshold function for a graph
property should have two-folds. They are corresponding to the
following two lemmas.
\begin{lem}\label{lem1} $rc_k(G(m,n,p))\geq d+2$ almost surely holds for every $p\leq ((\ln 2)^{1/d})/2\cdot p_1$.
\end{lem}
\pf Let $c_1=((\ln 2)^{1/d})/2$. Obviously $c_1<1$. Since
$$p\leq c_1 p_1=(1/2)\cdot\left((\ln (mn))^{1/d}/(m^{(d-1)/(2d)}n^{(d-1)/(2d)})\right),$$
by Theorem \ref{th2}, diam$(G(m,n,p))\geq d+2$ almost surely holds.
By
$$\mathbf{Pr}[rc_k(G(m,n,p))\geq d+2]\geq\mathbf{Pr}[\text{diam}(G(m,n,p))\geq d+2],$$
we get that for every $p\leq ((\ln 2)^{1/d})/2\cdot p_1$, $rc_k(G(m,n,p))\geq d+2$ almost surely holds.
\qed

\begin{lem}\label{lem2} $rc_k(G(m,n,p))\leq d+1$ almost surely holds for every $p\geq C_1\cdot p_1$.
\end{lem}
\pf Let $S=\{1,2,\ldots,d,d+1\}$ be a set of $d+1$ distinct colors. Randomly color the edges of
$G(m,n,p)$ with colors from $S$. By Theorem \ref{th3}, for every two distinct vertices
$u,v\in U(\text{or}\,\, u,v\in V)$ there are at least $2^{10d}c_0\log n$ internally vertex-disjoint
$uv$-paths of length exactly $d+1$. Let $P_1$ be such a $uv$-path. Under the random coloring,
the probability that $P_1$ is a rainbow path is
$$q_1=(d+1)!/(d+1)^{d+1}\geq ((d+1)/e)^{d+1}/(d+1)^{d+1}\geq 8^{-d},$$
by Stirling's formula. Meanwhile, for every $u\in U$ and $v\in V$ there are also at least
$2^{10d}c_0\log n$ internally vertex-disjoint $uv$-paths of length exactly $d$.
Let $P_2$ be such a $uv$-path. The probability that $P_2$ is a rainbow path is
$$q_2=d!/d^d\geq (d/e)^d/d^d\geq 4^{-d}.$$

Let $q=\min\{q_1,q_2\}\geq 8^{-d}$. Fix $u,v\in U(\text{or }\,
u,v\in V \,\text{or}\,\, u\in U,v\in V)$, we can estimate the upper
bound of the probability that there are at most $k-1$ such
$uv$-paths that are rainbow ones by
\begin{align*}
&{2^{10d}c_0\log n \choose k-1}(1-q)^{2^{10d}c_0\log n-(k-1)}\\
\leq &{2^{10d}c_0\log n \choose c_0\log n}(1-8^{-d})^{(2^{10d}-1)c_0\log n}\\
\leq &\left(\frac{2^{10d}c_0\log n\cdot e}{c_0\log n}\right)^{c_0\log n}\cdot2^{-8^{-d}(2^{10d}-1)c_0\log n}\\
=&(2^{10d}e)^{c_0\log n}\cdot 2^{-8^{-d}(2^{10d}-1)c_0\log n}\\
=&\left(\frac{e\cdot 2^{(10d+8^{-d})}}{2^{8^{-d}2^{10d}}}\right)^{c_0\log n}\\
\leq &n^{-100},
\end{align*}
where we apply the inequality ${n \choose k}\leq (\frac{ne}{k})^k$
\footnote[2]{We find that if in \cite{HL} He and Liang use this
inequality, their proof could be simplified significantly.}

By the Union Bound, with probability at least
$$1-\left[{m \choose
2}+{n \choose 2}+{m \choose 1}{n \choose 1}\right]n^{-100}\geq
1-n^{-90},$$ every two distinct vertices of $G(m,n,p)$ have at least
$k$ internally vertex-disjoint rainbow paths connecting them. This
implies that with probability at least $1-n^{-90}$, the event
$rc_k(G(m,n,p))\leq d+1$ happens, which gives precisely we want.
\qed

By the two lemmas above, it follows that
\begin{equation*}
p_1=(\log (mn))^{1/d}/(m^{(d-1)/(2d)}n^{(d-1)/(2d)})
\end{equation*}
is a sharp threshold function for the property $rc_k(G(m,n,p))\leq d+1$,
where $p_1n\geq p_1m \geq (\log n)^4$.

\noindent\textbf{Part 2: $d$ is even.}

Recall that $p_2=(\log n)^{1/d}/(m^{1/2}n^{(d-2)/(2d)})$. We
consider the random bipartite graph $G(m,n,p)$, where $m$ and $n$
satisfy that $p_2n^{1-\epsilon}\geq p_2m^{1-\epsilon} \geq (\log
n)^4$. The following two lemmas imply that $p_2$ is a sharp
threshold function for the property $rc_k(G(m,n,p))\leq d+1$.
\begin{lem}\label{lem13} $rc_k(G(m,n,p))\geq d+2$ almost surely holds for every $p\leq p_2$.
\end{lem}
\pf Let $c_2=1/((2\ln 2)^{1/d})$. Obviously, $c_2<1$. Since
\begin{eqnarray*}
p\leq p_2&=&(\log n)^{1/d}/(m^{1/2}n^{(d-2)/(2d)})\\
&=&\left(1/(2\ln2)^{1/d}\right)\cdot\left((2\ln n)^{1/d}/(m^{1/2}n^{(d-2)/(2d)})\right)\\
&=&c_2\cdot \left((2\ln n)^{1/d}/(m^{1/2}n^{(d-2)/(2d)})\right)
\end{eqnarray*}
by Theorem \ref{th2}, diam$(G(m,n,p))\geq d+2$ almost surely holds.
Then it follows that for every $p\leq p_2$, $rc_k(G(m,n,p))\geq d+2$
almost surely holds. \qed

Similar to the proof of Lemma \ref{lem2}, we can easily get the
following result.
\begin{lem}\label{lem14} $rc_k(G(m,n,p))\leq d+1$ almost surely holds for every $p\geq C_2\cdot p_2$.
\end{lem}

By the two lemmas above, we can conclude that
\begin{equation*}
p_2=(\log n)^{1/d}/(m^{1/2}n^{(d-2)/(2d)})
\end{equation*}
is a sharp threshold function for the property $rc_k(G(m,n,p))\leq
d+1$, where $m$ and $n$ satisfy that $p_2n^{1-\epsilon}
p_2m^{1-\epsilon} \geq (\log n)^4$.

Combining the two parts discussed above, we complete the proof of
Theorem \ref{th1}.

\section{The number of internally vertex-disjoint paths}

In this section, we prove Theorem \ref{th3}. We also divide our
proof into two parts according to the parity of $d$. We first give a
definition. An \emph{$(s,t)$-ary tree} with a designated root is a
tree such that every non-leaf vertex of even level has exactly $s$
children and every non-leaf vertex of odd level has exactly $t$
children, where we assume that the root is in zero-level. Obviously,
an $(s,t)$-ary tree and a $(t,s)$-ary tree of the same depth are
usually different trees.

\noindent\textbf{Part 1: $d$ is odd.}

Let $p=C_1p_1=C_1(\log (mn))^{1/d}/(m^{(d-1)/(2d)}n^{(d-1)/(2d)})$.
For every $u \in U$ and $S\subseteq V$ (or $u \in V$ and $S\subseteq
U$), let $X$ be the random variable counting the number of neighbors
of $u$ inside $S$.
\begin{lem}\label{lem3}
For every fixed $u,S$ such that $u\in U$, $S\subseteq V$ and
$|S|\geq n/2$ for sufficiently large $n$,
\begin{equation*}
\mathbf{Pr}[X\geq pn/10]\geq 1-2^{-\Omega(n^{1/d})}.
\end{equation*}
\end{lem}
\pf Denote by $S'$ any subset of $S$ with cardinality $n/2$. Let
$X_1$ be the random variable counting the number of neighbors of $u$
inside $S'$. Obviously, $X_1$ can be expressed as the sum of $n/2$
independent random variables, each of which taking $1$ with
probability $p$ and $0$ with probability $1-p$. Thus
$\mathbf{E}[X_1]=pn/2$. By the Chernoff-Hoeffding Bound, we have
\begin{equation*}
\mathbf{Pr}[X_1<(1-4/5)pn/2]\leq \exp\left(-(1/2)(4/5)^2(pn/2)\right)=2^{-\Omega(n^{1/d})},
\end{equation*}
By $X\geq X_1$, the event $X\geq pn/10$ happens with probability at least $1-2^{-\Omega(n^{1/d})}$,
which is precisely what we want.
\qed
\begin{lem}\label{lem4}
For every fixed $u,S$ such that $u\in V$, $S\subseteq U$ and
$|S|\geq m/2$ for sufficiently large $n$,
\begin{equation*}
\mathbf{Pr}[X\geq pm/10]\geq 1-n^{-\Omega(\log^3n)}.
\end{equation*}
\end{lem}

The proof is similar to that of Lemma \ref{lem3}. From Lemmas
\ref{lem3} and \ref{lem4}, it follows that $\mathbf{Pr}[X\geq
pn/\log n]\geq 1-2^{-\Omega(n^{1/d})}$ and $\mathbf{Pr}[X\geq
pm/\log m]\geq 1-n^{-\Omega(\log^3n)}$.
\begin{lem}\label{lem5} With probability at least $1-n^{-\Omega(1)}$, every two distinct
vertices of $U$ are connected by at least $2^{10d}c_0\log n$ internally vertex-disjoint
paths of length exactly $d+1$.
\end{lem}
\pf Fix $u,v\in U, u\neq v$. Consider the following process to
generate a $(pn/\log n,pm/\log m)$-ary tree of depth $d$ rooted at
$u$:
\begin{enumerate}[\indent Step 1.]
\item Let $T_0=\{u\}, i \leftarrow 1$, and $T_i\leftarrow \emptyset$.
\item If $i$ is odd, for every vertex $w\in T_{i-1}$, choose $pn/\log n$ distinct neighbors of $w$ from the
set $V\setminus(\cup^{i}_{j=0}T_j)$, and add them to $T_i$. (Note that $T_{i-1}\subseteq U$, $T_i$ is updated
every time after the processing of a vertex $w$, and in fact only when $j$ is odd $T_j\subseteq V$.)\\
If $i$ is even, for every vertex $w\in T_{i-1}$, choose $pm/\log m$
distinct neighbors of $w$ from the set $U\setminus(\{v\}\bigcup
\cup^i_{j=0}T_j)$, and add them to $T_i$.
\item Let $i\leftarrow i+1$. If $i\leq d$ then go to Step 2, otherwise stop.
\end{enumerate}

Of course, the process may fail during Step 2, since with nonzero
probability no neighbor of $w$ can be chosen as a candidate.
However, we will show that with high probability the tree can be
successfully constructed. Observe that when $j$ is even,
$T_j\subseteq U$, and when $j$ is odd, $T_j\subseteq V$. Thus,
$T_{d-1}\subseteq U$, $T_d\subseteq V$ and $|T_{d-1}|=(pm/\log
m)^{(d-1)/2}(pn/\log n)^{(d-1)/2}$, $|T_d|=(pm/\log
m)^{(d-1)/2}(pn/\log n)^{(d+1)/2}$. At any time during the process,
\begin{eqnarray*}
|\{v\}\cup \bigcup^i_{j=0,j\,\, \mbox{\footnotesize is even}}T_j|&\leq& 1+\sum^{d-1}_{j=0}|T_j|\leq d\cdot|T_{d-1}|\\
&=& d\cdot (pm/\log m)^{(d-1)/2}(pn/\log n)^{(d-1)/2}\\
&=&\frac{d\cdot C_1^{d-1}(\log (mn))^{(d-1)/d}}{(\log m)^{(d-1)/2}(\log n)^{(d-1)/2}}\cdot m^{(d-1)/2}n^{(d-1)/2}\\
&\leq& m/2
\end{eqnarray*}
and
\begin{eqnarray*}
|\bigcup^i_{j=0,j\,\, \mbox{\footnotesize is odd}}T_j|&\leq & 1+\sum_{j=0}^{d}{T_j}\leq (d+1)|T_d|\\
&=&(d+1)\cdot (pm/\log m)^{(d-1)/2}(pn/\log n)^{(d+1)/2}\\
&=&\frac{(d+1)\cdot C_1^d\log (mn)}{(\log m)^{(d-1)/2}(\log n)^{(d+1)/2}}\cdot n\\
&\leq& n/2
\end{eqnarray*}
for all sufficiently large $n$.

By Lemmas \ref{lem3} and \ref{lem4}, every execution of Step 2 fails
with probability at most $n^{-\Omega(\ln^4n)}$. Since Step 2 can be
executed for at most $(d+1)(pm/\log m)^{(d-1)/2}(pn/\log
n)^{(d+1)/2}$, we obtain that, with probability at least
\begin{equation*}
1-(d+1)\cdot(pm/\log m)^{(d-1)/2}(pn/\log n)^{(d+1)/2}\cdot n^{-\Omega(\log^3n)}=1-n^{-\Omega(\log^3n)},
\end{equation*}
the process can be successfully terminated.

Now we assume that $T$ has been successfully constructed. The number
of leaves in $T$ is exactly $|T_d|$. Let $Y$ be the random variable
counting the number of neighbors of $v$ inside $T_d$. It is obvious
that
\begin{equation*}
\mathbf{E}[Y]=p\cdot|T_d|=\frac{C_1^{d+1}(\log (mn))^{(d+1)/d}}{(\log m)^{(d-1)/2}(\log n)^{(d+1)/2}}
\cdot \frac{n^{(d+1)/(2d)}}{m^{(d-1)/(2d)}}\geq 10\cdot n^{1/(2d)}.
\end{equation*}
By the Chernoff-Hoeffding Bound, we get
\begin{equation*}
\mathbf{Pr}[Y<n^{1/(2d)}]\leq \exp(-(1/2)(9/10)^2\cdot10n^{1/(2d)})\leq 2^{-n^{1/(4d)}}
\end{equation*}

For every $w\in T_1$, define as the vice-tree $T_w$ of $T$ the
subtree of $T$ of depth $d-1$ rooted at $w$. Notice that every
vice-tree contains $(pm/\log m)^{(d-1)/2}(pn/\log n)^{(d-1)/2}$
leaves. For each vice-tree $T_w$, let $Z_w$ be the random variable
counting the number of neighbors of $v$ inside the set of leaves of
$T_w$. Then we have
\begin{eqnarray*}
\mathbf{Pr}[Z_w\geq n^{1/(10d)}]&\leq& {(pm/\log m)^{(d-1)/2}(pn/\log n)^{(d-1)/2} \choose n^{1/(10d)}}\cdot p^{n^{1/(10d)}}\\
&\leq& \left(\frac{(pm/\log m)^{(d-1)/2}(pn/\log n)^{(d-1)/2}\cdot e}{n^{1/(10d)}}\right)^{n^{1/(10d)}}\cdot p^{n^{1/(10d)}}\\
&=& \left(\frac{C_1^d\log (mn)\cdot e}{(\log m)^{(d-1)/2}(\log n)^{(d-1)/2}\cdot n^{1/(10d)}}\right)^{n^{1/(10d)}} \\
&\leq& n^{-100},
\end{eqnarray*}
where we apply the inequality ${n \choose k}\leq (\frac{ne}{k})^k$. By applying the Union Bound, we get
\begin{equation*}
\mathbf{Pr}[\bigvee_{w\in T_1}Z_w]\leq (pn/\log n)\cdot n^{-100}\leq n^{-90}.
\end{equation*}

Combined with previous estimations, we derive that with probability
at least
\begin{equation*}
1-n^{-\Omega(\log^3n)}-2^{-n^{1/(4d)}}-n^{-90}\geq 1-n^{-80},
\end{equation*}
the following three events simultaneously happen:
\begin{enumerate}[\indent 1.]
\item the tree $T$ is successfully constructed,
\item $v$ has at least $n^{1/(2d)}$ neighbors inside the set of leaves of
$T$,
\item every vice-tree $T_w$ contains at most $n^{1/(10d)}$ leaves that are neighbors of $v$.
\end{enumerate}

It is clear that each neighbor $v'$ of $v$ inside $T_d$ induces a
$uv$-path of length $d+1$. If two neighbors $v'$ and $v''$ of $v$
belong to distinct vice-trees, then the corresponding two $uv$-paths
are internally vertex-disjoint. When all these three events happen,
we can choose $n^{1/(2d)}/n^{1/(10d)}=n^{2/(5d)}\geq 2^{10d}c_0\log
n$ neighbors of $v$ inside $T_d$, every two of which are from
different vice-trees. Thus we can immediately obtain at least
$2^{10d}c_0\log n$ internally vertex-disjoint $uv$-paths.

By using the Union Bound again, it then follows that, with
probability at least
$$1-{m \choose 2}\cdot n^{-80}=1-n^{-\Omega(1)},$$
every two distinct vertices of $U$ are
connected by at least $2^{10d}c_0\log n$ internally vertex-disjoint
paths of length exactly $d+1$. This completes the proof of Lemma
\ref{lem5}. \qed

\begin{lem}\label{lem6} With probability at least $1-n^{-\Omega(1)}$, every two distinct
vertices of $V$ are connected by at least $2^{10d}c_0\log n$ internally vertex-disjoint
paths of length exactly $d+1$.
\end{lem}
\pf The proof is similar to that of Lemma \ref{lem5}. Here we only
point out the differences. Fix $u,v\in V, u\neq v$. We first
construct a $\left(pm/((\log m)^{2/(d-1)}),pn/((\log
n)^{2/(d-1)})\right)$-ary tree of depth $d$ rooted at $u$. We can
easily determinate that with probability at least
$1-n^{-\Omega(\log^3n)}$ the tree $T$ can be successfully
constructed. $T_d$ is just the set of leaves of $T$ and
$$T_d=\left(pm/((\log m)^{2/(d-1)})\right)^{(d+1)/2}\cdot
\left(pn/((\log n)^{2/(d-1)})\right)^{(d-1)/2}.$$

Let $Y$ be the random variable counting the number of neighbors $v$
inside $T_d$. It is obvious that
\begin{equation*}
\mathbf{E}[Y]=p\cdot|T_d|=\frac{C_1^{d+1}(\log (mn))^{(d+1)/d}}{(\log m)^{(d+1)/(d-1)}\log n}
\cdot \frac{m^{(d+1)/(2d)}}{n^{(d-1)/(2d)}}\geq C_1^2(\log n)^{2}.
\end{equation*}
As before, we have
\begin{eqnarray*}
\mathbf{Pr}\left[Y<C_1(\log n)^2\right]&\leq& \mathbf{Pr}\left[Y<(1-(C_1-1)/C_1)\cdot C_1(\log n)^2\right]\\
&\leq& \exp\left(-(1/2)((C_1-1)/C_1)^2\cdot C_1(\log n)^2\right)\\
&\leq& n^{-O(\log n)}
\end{eqnarray*}

Notice that every vice-tree contains $\left(pm/((\log
m)^{2/(d-1)})\right)^{(d-1)/2}\cdot \left(pn/((\log
n)^{2/(d-1)})\right)^{(d-1)/2}$
$=\left(p^{d-1}m^{(d-1)/2}n^{(d-1)/2}\right)/\left(\log m\cdot \log
n\right)$ leaves. For each vice-tree $T_w$, let $Z_w$ be the random
variable counting the number of neighbors $v$ inside the set of
leaves of $T_w$. Then we have
\begin{eqnarray*}
\mathbf{Pr}[Z_w\geq \log n] &\leq& {\left(p^{d-1}m^{(d-1)/2}n^{(d-1)/2}\right)/\left(\log m
\cdot \log n\right) \choose \log n}\cdot p^{\log n}\\
&\leq& \left(\frac{p^dm^{(d-1)/2}n^{(d-1)/2}e}{\log m\cdot \log^2 n}\right)^{\log n}\\
&=& \left(\frac{C_1^d\log (mn)\cdot e}{\log m\cdot\log^2 n}\right)^{\log n} \\
&\leq& n^{-O(\log\log n)},
\end{eqnarray*}
and
\begin{equation*}
\mathbf{Pr}[\bigvee_{w\in T_1}Z_w]\leq (pm/(\log m)^{2/(d-1)})\cdot n^{-O(\log\log n)}=n^{-O(\log\log n)}.
\end{equation*}

Since $C_1\log^2n/\log n=2^{10d}c_0\log n$, combined with the estimations above, we derive that with
probability at least
$$1-n^{-\Omega(\log^3n)}-n^{-O(\log n)}-n^{-O(\log\log n)}= 1-n^{-O(\log\log n)},$$
there are at least $2^{10d}c_0\log n$ internally vertex-disjoint
$uv$-paths.

Therefore, we can easily obtain that with probability at least
$$1-{n \choose 2}\cdot n^{-O(\log\log n)}=1-n^{-O(\log\log n)},$$
every two distinct vertices of $V$ have at least $2^{10d}c_0\log n$
internally vertex-disjoint paths of length $d+1$ connecting them.
\qed

\begin{lem}\label{lem7} With probability at least $1-n^{-\Omega(1)}$, every two distinct vertices of
different partites are connected by at least $2^{10d}c_0\log n$
internally vertex-disjoint paths of length exactly $d$.
\end{lem}
\pf Similarly, fix $u\in U, v\in V$. We first construct a
$(pn/10,pm/10)$-ary tree of depth $d-1$ rooted at $u$. It can also
be estimated that with probability at least $1-n^{-\Omega(\log^3n)}$
the tree $T$ can be successfully constructed.

Let $Y$ be the random variable counting the number of neighbors of
$v$ inside $T_{d-1}$ which is just the set of leaves of $T$. Then
\begin{eqnarray*}
&|T_{d-1}|=(pm/10)^{(d-1)/2}(pn/10)^{(d-1)/2},\\
&\mathbf{E}[Y]=p\cdot|T_{d-1}|= 10\cdot (C_1/10)^d\cdot \log (mn),\\
&\mathbf{Pr}\left[Y<(C_1/10)^d\log (mn)\right]\leq n^{-10}.
\end{eqnarray*}

Notice that every vice-tree contains
$(pm/10)^{(d-1)/2}(pn/10)^{(d-3)/2}$ leaves. For each vice-tree
$T_w$, let $Z_w$ be the random variable counting the number of
neighbors $v$ inside the set of leaves of $T_w$. Then we have
\begin{equation*}
\mathbf{Pr}[Z_w\geq 10d] \leq {(pm/10)^{(d-1)/2}\cdot (pn/10)^{(d-3)/2} \choose 10d}\cdot p^{10d}\leq O(n^{-5}),
\end{equation*}
\begin{equation*}
\mathbf{Pr}[\bigvee_{w\in T_1}Z_w]\leq (pn/10)\cdot O(n^{-5})\leq O(n^{-4}).
\end{equation*}

Since $((C_1/10)^d/10d)\log (mn)\geq 2^{10d}c_0\log n$, combined with the estimations above,
we derive that with probability at least
$$1-n^{-\Omega(\log^3n)}-n^{-10}-O(n^{-4})\geq 1-O(n^{-3}),$$
there are at least $2^{10d}c_0\log n$ internally vertex-disjoint
$uv$-paths.

It follows that with probability at least $1-mn\cdot O(n^{-3})\geq
1-n^{-1}$ every two distinct vertices of different partites have at
least $2^{10d}c_0\log n$ internally vertex-disjoint paths of length
$d$ connecting them.
\qed

Now we have seen that Theorem \ref{th3} is true for the case that
$d$ is odd.

\noindent\textbf{Part 2: $d$ is even.}

Let $p=C_2p_2=C_2(\log n)^{1/d}/(m^{1/2}n^{(d-2)/(2d)})$. For every
$u \in U$ and $S\subseteq V$ (or $u \in V$ and $S\subseteq U$), let
$X$ be the random variable counting the number of neighbors of $u$
inside $S$. We have the following results similar to Lemmas
\ref{lem3} and \ref{lem4}.
\begin{lem}\label{lem8}
For every fixed $u,S$ such that $u\in U$, $S\subseteq V$ and $|S|\geq n/2$ for sufficiently large $n$,\\
\begin{equation*}
\mathbf{Pr}[X\geq pn/10]\geq 1-2^{-\Omega(n^{1/d})}.
\end{equation*}
\end{lem}
\begin{lem}\label{lem9}
For every fixed $u,S$ such that $u\in V$, $S\subseteq U$ and $|S|\geq m/2$ for sufficiently large $n$,\\
\begin{equation*}
\mathbf{Pr}[X\geq pm/10]\geq 1-n^{-\Omega(\log^3n)}.
\end{equation*}
\end{lem}

Lemmas \ref{lem8} and \ref{lem9} also imply that $\mathbf{Pr}[X\geq
pn/\log n]\geq 1-2^{-\Omega(n^{1/d})}$ and $\mathbf{Pr}[X\geq
pm/\log m]\geq 1-n^{-\Omega(\log^3n)}$. The proofs of the following
three lemmas are similar to those of Lemmas \ref{lem5}, \ref{lem6}
and \ref{lem7}, but the estimations are different. So we only sketch
the proofs of them and list the results of their estimations.
\begin{lem}\label{lem10} With probability at least $1-n^{-\Omega(1)}$,
every two distinct vertices of $U$ are connected by at least $2^{10d}c_0\log n$
internally vertex-disjoint paths of length exactly $d$.
\end{lem}
\pf Fix $u,v\in U, u\neq v$. We construct a $(pn/10,pm/10)$-ary tree
of depth $d-1$ rooted at $u$. Then we can similarly give the
estimations. We obtain that with probability at least
$$1-n^{-\Omega(\log^3n)}-n^{-10}-n^{-5}\geq 1-n^{-4},$$ the
following three events simultaneously happen:
\begin{enumerate}[\indent 1.]
\item the tree $T$ is successfully constructed,
\item $v$ has at least $((C_2/10)^d\log n)\cdot (n/m)$ neighbors inside the set of leaves of
$T$,
\item every vice-tree $T_w$ contains at most $10d\cdot (n/m)$ leaves that are neighbors of $v$.
\end{enumerate}

Since
$$\frac{((C_2/10)^d\log n)\cdot (n/m)}{10d\cdot (n/m)}\geq 2^{10d}c_0\mbox{log}n,$$
and
$$1-{n \choose 2}\cdot n^{-4}=1-n^{-\Omega(1)},$$
every two distinct vertices of $U$ are connected
by at least $2^{10d}c_0\log n$ internally vertex-disjoint paths of
length exactly $d$.
\qed
\begin{lem}\label{lem11} With probability at least $1-n^{-\Omega(1)}$, every two distinct vertices of
$V$ are connected by at least $2^{10d}c_0\log n$ internally vertex-disjoint paths of length exactly $d$.
\end{lem}
\pf Fix $u,v\in V, u\neq v$. In this case, we construct a
$(pm/10,pn/10)$-ary tree of depth $d-1$ rooted at $u$. We can
determine that with probability at least
$$1-n^{-\Omega(\log^3n)}-n^{-10}-O(n^{-8})\geq 1-O(n^{-7}),$$ the
following three events simultaneously happen:
\begin{enumerate}[\indent 1.]
\item the tree $T$ is successfully constructed,
\item $v$ has at least $(C_2/10)^d\log n$ neighbors inside the set of leaves of
$T$,
\item every vice-tree $T_w$ contains at most $(10d)/\epsilon$ leaves that are neighbors of $v$.
\end{enumerate}

Since
$$\frac{(C_2/10)^d\log n}{(10d)/\epsilon}\geq 2^{10d}c_0\log n,$$
and
$$1-{n \choose 2}\cdot O(n^{-7})\geq 1-n^{-\Omega(1)},$$
every two distinct vertices of $V$ are connected by at least
$2^{10d}c_0\mbox{log}\,n$ internally vertex-disjoint paths of length
exactly $d$.
\qed
\begin{lem}\label{lem12} With probability at least $1-n^{-\Omega(1)}$,
every two distinct vertices of different partites are connected by
at least $2^{10d}c_0\log n$ internally vertex-disjoint paths of
length exactly $d+1$.
\end{lem}
\pf Fix $u\in U, v\in V$. In this case, the estimations are more
complicated than the previous cases. We construct a $(pm/\log
m,pn/\log n)$-ary tree of depth $d$ rooted at $v$. We can determine
that with probability at least
$$1-n^{-\Omega(\log^3n)}-n^{-10}-O(n^{-8})\geq 1-O(n^{-6}),$$
the following three events simultaneously happen:
\begin{enumerate}[\indent 1.]
\item the tree $T$ is successfully constructed,
\item $u$ has at least
$$L_1=\frac{C_2^d(\log n)^{(d+1)/d}}{(\log m)^{d/2}(\log n)^{d/2}}\cdot \frac{n^{1/2+1/d}}{m^{1/2}}$$
neighbors inside the set of leaves of $T$,
\item every vice-tree $T_w$ contains at most
$$L_2=\frac{C_2^d\log n}{(\log m)^{d/2-1}(\log n)^{d/2}}\cdot \frac{n}{m^{1-\epsilon/2}}$$
leaves that are neighbors of $v$.
\end{enumerate}

Since
$$L_1/L_2\geq 2^{10d}c_0\log n,$$
and
$$1-{n \choose 2}\cdot O(n^{-6})\geq 1-n^{-\Omega(1)},$$
every two distinct vertices of different partites are connected by
at least $2^{10d}c_0\log n$ internally vertex-disjoint paths of
length exactly $d$. \qed

Now we can see that Theorem \ref{th3} is also true for the case that
$d$ is even.

Combining the two parts discussed above, the proof of Theorem
\ref{th3} is thus completed.

\end{document}